\documentstyle[bezier]{article}

\def\vpint{\int_{-\infty}^\infty\hspace{-20.5pt}-\hspace{11.5pt}}  
\catcode`\@=11
\ifcase\@ptsize
 \font\tenmsa=msam10
 \font\sevenmsa=msam7
 \font\fivemsa=msam5
 \font\tenmsb=msbm10
 \font\sevenmsb=msbm7
 \font\fivemsb=msbm5
 \font\teneu=eufm10
 \font\seveneu=eufm7
 \font\fiveeu=eufm5
 \font\tenib=cmmib10
 \font\sevenib=cmmib7
 \font\fiveib=cmmib5
\or
 \font\tenmsa=msam10 scaled \magstephalf
 \font\sevenmsa=msam7 scaled \magstephalf
 \font\fivemsa=msam5 scaled \magstephalf
 \font\tenmsb=msbm10 scaled \magstephalf
 \font\sevenmsb=msbm7 scaled \magstephalf
 \font\fivemsb=msbm5  scaled \magstephalf
 \font\teneu=eufm10  scaled \magstephalf
 \font\seveneu=eufm7  scaled \magstephalf
 \font\fiveeu=eufm5   scaled \magstephalf
 \font\tenib=cmmib10  scaled \magstephalf
 \font\sevenib=cmmib7  scaled \magstephalf
 \font\fiveib=cmmib5   scaled \magstephalf
\or
 \font\tenmsa=msam10 scaled \magstep1
 \font\sevenmsa=msam7 scaled \magstep1
 \font\fivemsa=msam5  scaled \magstep1
 \font\tenmsb=msbm10 scaled \magstep1
 \font\sevenmsb=msbm7 scaled \magstep1
 \font\fivemsb=msbm5  scaled \magstep1
 \font\teneu=eufm10   scaled \magstep1
 \font\seveneu=eufm7 scaled \magstep1
 \font\fiveeu=eufm5 scaled \magstep1
 \font\tenib=cmmib10     scaled \magstep1
 \font\sevenib=cmmib7   scaled \magstep1
 \font\fiveib=cmmib5   scaled \magstep1
\fi
\newfam\msafam
\textfont\msafam=\tenmsa  \scriptfont\msafam=\sevenmsa
  \scriptscriptfont\msafam=\fivemsa
\newfam\msbfam
\textfont\msbfam=\tenmsb  \scriptfont\msbfam=\sevenmsb
  \scriptscriptfont\msbfam=\fivemsb
\def\Bbb{\ifmmode\let\next\Bbb@\else
 \def\next{\errmessage{Use \string\Bbb\space only in math mode}}\fi\next}
\def\Bbb@#1{{\Bbb@@{#1}}}
\def\Bbb@@#1{\fam\msbfam#1}
\newfam\eufam
\textfont\eufam=\teneu  \scriptfont\eufam=\seveneu
  \scriptscriptfont\eufam=\fiveeu
\def\frak{\ifmmode\let\next\frak@\else
 \def\next{\errmessage{Use \string\frak\space only in math mode}}\fi\next}
\def\frak@#1{{\frak@@{#1}}}
\def\frak@@#1{\fam\eufam#1}
\newfam\ibfam
\textfont\ibfam=\tenib  \scriptfont\ibfam=\sevenib
  \scriptscriptfont\ibfam=\fiveib
\def\bold{\ifmmode\let\next\bold@\else
 \def\next{\errmessage{Use \string\bold\space only in math mode}}\fi\next}
\def\bold@#1{{\bold@@{#1}}}
\def\bold@@#1{\fam\ibfam#1}
\def\hexnumber@#1{\ifcase#1 0\or 1\or 2\or 3\or 4\or 5\or 6\or 7\or 8\or
 9\or A\or B\or C\or D\or E\or F\fi}
\def\newsymbolb#1#2#3#4{\mathchardef#1="#2\hexnumber@\msbfam#3#4}
\def\newsymbola#1#2#3#4{\mathchardef#1="#2\hexnumber@\msafam#3#4}
\font\fraksect=eufm10 scaled 1728

\font\frakssect=eufm10 scaled 1440
\def\hybrid{\topmargin 0pt      \oddsidemargin 0pt
        \headheight 0pt \headsep 0pt
        \textwidth 160true mm       
        \textheight 231true mm         
        \marginparwidth 0.0in
        \parskip 0pt plus 1pt   \jot = 1.5ex}
\catcode`\@=11
\def\marginnote#1{}

\newcount\hour
\newcount\minute
\newtoks\amorpm
\hour=\time\divide\hour by60
\minute=\time{\multiply\hour by60 \global\advance\minute by-\hour}
\edef\standardtime{{\ifnum\hour<12 \global\amorpm={am}%
        \else\global\amorpm={pm}\advance\hour by-12 \fi
        \ifnum\hour=0 \hour=12 \fi
        \number\hour:\ifnum\minute<10 0\fi\number\minute\the\amorpm}}
\edef\militarytime{\number\hour:\ifnum\minute<10 0\fi\number\minute}

\def\draftlabel#1{{\@bsphack\if@filesw {\let\thepage\relax
   \xdef\@gtempa{\write\@auxout{\string
      \newlabel{#1}{{\@currentlabel}{\thepage}}}}}\@gtempa
   \if@nobreak \ifvmode\nobreak\fi\fi\fi\@esphack}
        \gdef\@eqnlabel{#1}}
\def\@eqnlabel{}
\def\@vacuum{}
\def\draftmarginnote#1{\marginpar{\raggedright\scriptsize\tt#1}}

\def\draft{\oddsidemargin -.5truein
        \def\@oddfoot{\sl preliminary draft \hfil
        \rm\thepage\hfil\sl\today\quad\militarytime}
        \let\@evenfoot\@oddfoot \overfullrule 3pt
        \let\label=\draftlabel
        \let\marginnote=\draftmarginnote
   \def\@eqnnum{(\theequation)\rlap{\kern\marginparsep\tt\@eqnlabel}%
\global\let\@eqnlabel\@vacuum}  }
\newcounter{app}
\newcounter{sapp}[app]

\newcommand{\app}[1]{
\refstepcounter{app}{\vspace{7mm}
\noindent\Large\bf Appendix
\theapp.
 \ #1 \par \vspace{5mm}}
\setcounter{equation}{0}
\def\theequation{\Alph{app}.\arabic{equation}}}

\makeatletter
\newdimen\normalarrayskip              
\newdimen\minarrayskip                 
\normalarrayskip\baselineskip
\minarrayskip\jot
\newif\ifold             \oldtrue            
\def\arraymode{\ifold\relax\else\displaystyle\fi} 
\def\eqnumphantom{\phantom{(\theequation)}}     
\def\@arrayskip{\ifold\baselineskip\z@\lineskip\z@
     \else
     \baselineskip\minarrayskip\lineskip2\minarrayskip\fi}
\def\@arrayclassz{\ifcase \@lastchclass \@acolampacol \or
\@ampacol \or \or \or \@addamp \or
   \@acolampacol \or \@firstampfalse \@acol \fi
\edef\@preamble{\@preamble
  \ifcase \@chnum
     \hfil$\relax\arraymode\@sharp$\hfil
     \or $\relax\arraymode\@sharp$\hfil
     \or \hfil$\relax\arraymode\@sharp$\fi}}
\def\@array[#1]#2{\setbox\@arstrutbox=\hbox{\vrule
     height\arraystretch \ht\strutbox
     depth\arraystretch \dp\strutbox
     width\z@}\@mkpream{#2}\edef\@preamble{\halign \noexpand\@halignto
\bgroup \tabskip\z@ \@arstrut \@preamble \tabskip\z@ \cr}%
\let\@startpbox\@@startpbox \let\@endpbox\@@endpbox
  \if #1t\vtop \else \if#1b\vbox \else \vcenter \fi\fi
  \bgroup \let\par\relax
  \let\@sharp##\let\protect\relax
  \@arrayskip\@preamble}
\def\eqnarray{\stepcounter{equation}%
              \let\@currentlabel=\theequation
              \global\@eqnswtrue
              \global\@eqcnt\z@
              \tabskip\@centering
              \let\\=\@eqncr
              $$%
 \halign to \displaywidth\bgroup
    \eqnumphantom\@eqnsel\hskip\@centering
    $\displaystyle \tabskip\z@ {##}$%
    &\global\@eqcnt\@ne \hskip 2\arraycolsep
         $\displaystyle\arraymode{##}$\hfil
    &\global\@eqcnt\tw@ \hskip 2\arraycolsep
         $\displaystyle\tabskip\z@{##}$\hfil
         \tabskip\@centering
    &{##}\tabskip\z@\cr}
\makeatother
\relax
\hybrid
\begin{document}
\def\bea{\begin{eqnarray}}
\def\eea{\end{eqnarray}}
\def\beq{\begin{equation}}          \def\bn{\beq}
\def\eeq{\end{equation}}            \def\ed{\eeq}
\def\nn{\nonumber}                  \def\g{\gamma}
\def\Uq{U_q(\widehat{\frak{sl}}_2)}
\def\Uqp{U_q(\widehat{\frak{sl}}'_2)}
\def\Uqd{U^{*}_q(\widehat{\frak{sl}}_2)}
\def\uq{U_q({sl}_2)}
\def\uqd{U^*_q({sl}_2)}
\def\slaff{\frak{sl}^\prime_2}
\def\aff{\widehat{\frak{sl}}_2}
\def\ot{\otimes}
\def\sk#1{\left({#1}\right)}
\def\id{\mbox{\rm id}}
\def\tr{\mbox{\rm tr}}
\def\tah{\mbox{\rm th}}
\def\sh{\mbox{\rm sh}}
\def\ch{\mbox{\rm ch}}
\def\ctg{\mbox{\rm ctg}}
\def\cth{\mbox{\rm cth}}
\def\tg{\mbox{\rm tg}}
\def\th{\mbox{\rm th}}
\def\qdet{\mbox{\rm q-det}}
\def\Re{{\rm Re}\,}
\def\Im{{\rm Im}\,}
\def\sn{\mbox{\rm sn}}
\def\cn{\mbox{\rm cn}}
\def\dn{\mbox{\rm dn}}
\def\snh{\mbox{\rm snh}}
\def\cnh{\mbox{\rm cnh}}
\def\dnh{\mbox{\rm dnh}}
\def\RR{\Bbb{R}}
\def\ZZ{\Bbb{Z}}
\def\CC{\Bbb{C}}
\def\r#1{\mbox{(}\ref{#1}\mbox{)}}
\def\d{\delta}
\def\D{\Delta}
\def\da{{\partial_\alpha}}
\let\da=p
\def\Ps{\Psi^{*}}
\def\R{{\cal R}}
\def\Ga#1{\Gamma\left(#1\right)}
\def\si#1{\sin\left(#1\right)}
\def\ex#1{\exp\left(#1\right)}
\def\ep{\varepsilon}
\def\eps{\epsilon}
\def\ve{\ep}
\def\fract#1#2{{\mbox{\footnotesize $#1$}\over\mbox{\footnotesize $#2$}}}
\def\stackreb#1#2{\ \mathrel{\mathop{#1}\limits_{#2}}}
\def\stackreu#1#2{\ \mathrel{\mathop{#1}\limits^{#2}}}
\def\res#1{\stackreb{\mbox{\rm res}}{#1}}
\def\lim#1{\stackreb{\mbox{\rm lim}}{#1}}
\def\Res#1{\stackreb{\mbox{\rm Res}}{#1}}
\let\dis=\displaystyle
\def\ee{{\rm e}}
\def\D{\Delta}
\renewcommand{\theequation}{{\thesection}.{\arabic{equation}}}
\def\Y-{\widehat{Y}^-}
\font\fraksect=eufm10 scaled 1728
\font\frakssect=eufm10 scaled 1440
\def\DYsect{\widehat{DY(\hbox{\fraksect sl}_2)}}
\def\Ael{{\cal A}_{\tih,\eta}(\widehat{\frak{sl}_2})}
\def\Aelx{{\cal A}_{\tih,1/\xi}(\widehat{\frak{sl}_2})}
\def\Aelxsect{{\cal A}_{\tih,1/\xi}(\widehat{\hbox{\fraksect sl}_2})}
\def\Aelxssect{{\cal A}_{\tih,1/\xi}(\widehat{\hbox{\frakssect sl}_2})}
\def\aelxc{\frak{a}_{1/\xi,c}(\widehat{\frak{sl}_2})}
\def\aelc{\frak{a}_{\eta,c}(\widehat{\frak{sl}_2})}
\def\aelxcsect{\hbox{\fraksect a}_{1/\xi,c}(\widehat{\hbox{\fraksect sl}_2})}
\def\aelxcssect{\hbox{\frakssect a}_{1/\xi,c}
(\widehat{\hbox{\frakssect sl}_2})}
\def\aelcsect{\hbox{\fraksect a}_{\eta,c}(\widehat{\hbox{\fraksect sl}_2})}
\def\aelcssect{\hbox{\frakssect a}_{\eta,c}(\widehat{\hbox{\frakssect sl}_2})}
\def\aelcsseco{\hbox{\frakssect a}_{\eta}^0(\widehat{\hbox{\frakssect sl}_2})}
\def\Aelsect{{\cal A}_{\tih,\eta}(\widehat{\hbox{\fraksect sl}_2})}
\def\Aelssect{{\cal A}_{\tih,\eta}(\widehat{\hbox{\frakssect sl}_2})}
\def\ael{{\frak{a}}_{\eta}(\widehat{\frak{sl}_2})}
\def\aelo{{\frak{a}}_{\eta}^0(\widehat{\frak{sl}_2})}
\def\aelt{\tilde{\frak{a}}_{\eta}^0(\widehat{\frak{sl}_2})}
\def\aelp{{\frak{a}}^+_{\eta}(\widehat{\frak{sl}_2})}
\def\aelm{{\frak{a}}^-_{\eta}(\widehat{\frak{sl}_2})}
\def\aelpm{{\frak{a}}^\pm_{\eta}(\widehat{\frak{sl}_2})}
\def\aelmp{{\frak{a}}^\mp_{\eta}(\widehat{\frak{sl}_2})}
\def\ap{{\frak{a}}_{p}(\widehat{\frak{sl}_2})}
\def\apo{{\frak{a}}_{p}^0(\widehat{\frak{sl}_2})}
\def\app{{\frak{a}}^+_{p}(\widehat{\frak{sl}_2})}
\def\apm{{\frak{a}}^-_{p}(\widehat{\frak{sl}_2})}
\def\apssect{{\hbox{\frakssect a}}_{p}(\widehat{\hbox{\frakssect sl}_2})}
\def\aelx{{\frak{a}}_{1/\xi}(\widehat{\frak{sl}_2})}
\def\aelsect{\hbox{\fraksect{a}}_{\eta}(\widehat{\hbox{\fraksect sl}_2})}
\def\aelseco{\hbox{\fraksect{a}}_{\eta}^0(\widehat{\hbox{\fraksect sl}_2})}
\def\aelssect{\hbox{\frakssect{a}}_{\eta}(\widehat{\hbox{\frakssect sl}_2})}
\def\Apq{{\cal A}_{q,p}(\widehat{\frak{sl}_2})}
\def\Apqsect{{\cal A}_{q,p}(\widehat{\hbox{\fraksect sl}_2})}
\def\ap{{\frak{a}}_{p}(\widehat{\frak{sl}_2})}
\def\apsect{\hbox{\fraksect{a}}_{p}(\widehat{\hbox{\fraksect sl}_2})}
\def\DY{\widehat{DY(\frak{sl}_2)}}
\def\Yd{\DY}
\def\Ydd{\DY}
\let\z=z
\let\b=z
\def\u{{u}}
\def\v{{v}}
\def\g{\gamma}
\def\la{\lambda}
\let\hsp=\qquad
\def\he{{\hat e}}
\def\hf{{\hat f}}
\def\hh{{\hat h}}
\def\ha{{\hat a}}
\def\hb{{\hat b}}
\def\hk{{\hat t}}
\def\hkp{{\hat t}}
\def\hhh{{\hat h}}
\def\vvv{\overline{\varphi}}
\def\tle{{\tilde e}}
\def\tlf{{\tilde f}}
\def\tlh{{\tilde h}}
\def\intt{\int_{-\infty}^\infty}
\def\la{\lambda}
\def\tih{{\hbar}}
\def\FDY{F\left[\DY\right]}
\def\FDYsect{F\left[\DYsect\right]}
\def\stackupb#1#2#3{\ \mathrel{\mathop{#1}\limits_{#2}^{#3}}}
\def\feq#1#2{\stackupb{\ravnodots}{#1}{#2}}
\def\u{{u\over\xi}}
\def\H{{\cal H}}
\def\cint{\int_\infty^{0+}}
\def\nint{\int^{+\infty}_{0}}
\def\mint{\int_{-\infty}^{0}}
\def\rvac{|\mbox{vac}\rangle}
\def\lvac{\langle\mbox{vac}|}
\begin{center}
\hfill ITEP-TH-1/97\\
\hfill RIMS-1139\\
\hfill q-alg/9703043\\
\bigskip\bigskip
{\Large\bf Classical Limit of the Scaled  Elliptic Algebra $\Aelsect$}\\
\bigskip
\bigskip
{\large 
S. Khoroshkin$^*$\footnote{E-mail: khor@heron.itep.ru},
D. Lebedev$^*$\footnote{E-mail: dlebedev@vitep5.itep.ru},
S. Pakuliak$^{*\star\diamond}$\footnote{E-mail: pakuliak@thsun1.jinr.ru},
A. Stolin$^\circ$\footnote{E-mail: astolin@math.chalmers.se},
V. Tolstoy$^\dagger$\footnote{E-mail: tolstoy@anna19.npi.msu.su}}
\end{center}
\bigskip

$^*${\it Institute of Theoretical \& Experimental Physics,
117259 Moscow, Russia}

$^\star${\it Bogoliubov Laboratory of Theoretical Physics, JINR,
141980 Dubna, Moscow region, Russia}

$^\diamond${\it Bogoliubov Institute of Theoretical Physics, 252143
Kiev, Ukraine}

$^\circ${\it Department of Mathematics,
G\"oteborg University, S-41296 G\"oteborg, Sweden}

$^\dagger${\it Institute of Nuclear Physics,
Moscow State University, 119899 Moscow, Russia}


{\bf AMS Subj. Class.:} 17B65 (Primary) 17B37, 46F20 (Secondary)

\bigskip
\bigskip

\begin{abstract}
The classical limit of the scaled elliptic algebra $\Ael$ is investigated.
The limiting Lie algebra is described in two equivalent ways: as a central
 extension of the algebra of 
 generalized automorphic $sl_2$ valued functions on a strip and
 as an extended algebra of decreasing automorphic $sl_2$ valued functions on
 the real line. A bialgebra structure and an infinite-dimensional 
 representation
 in the Fock space are studied. The classical limit of elliptic algebra
 $\Apq$ is also briefly presented.

\end{abstract}

\setcounter{section}{0}
\setcounter{equation}{0}
\setcounter{footnote}{0}
\section{Introduction}

The elliptic  algebra $\Apq$ has been introduced in the paper
\cite{MJ}. Its definition was induced mostly by the bosonization
formulas for massive integrable field theories, proposed in \cite{L}.
Then it was shown in \cite{MJ2} that in the framework of the
representation theory of the algebra $\Ael$ which is a scaling
limit of the algbera $\Apq$ one can obtain integral representations
for the correlation functions in the $XXZ$-model in massless regime.

The structure of the algebra $\Ael$ is rather unusual.
 For its precise definition one should introduce a continuous family of
 generators being Fourier harmonics of the Gauss coordinates 
 of the $L$-operator.
 The elements of the algebra are  formal integrals over the
 generators with certain conditions on analyticity and on asymptotics
 of the coefficients. Next,
 both algebras $\Apq$ and $\Ael$ are not Hopf algebras but form
a Hopf family of  algebras (see section 4)  and even at level
$c=0$ when these algebras become usual Hopf algebras, they do not
have the structure of a double which can be reconstructed
either in the Yangian limit $\eta\to 0$ from $\Ael$ or in the quantum
 affine limit $(p\to 0)$ of $\Apq$.
 We are interested in corresponding properties of the limiting
 (with respect to $\Ael$) classical
 algebra $\ael$ which we consider in details. We observe also the elliptic
 case and the rational degeneration of Lie algebra $\ael$.

 The limiting algebra $\ael$ can be
 described in two ways. First, it can be
realized as
 a central extension of
the algebra of automorphic $sl_2$-valued
 generalized functions on a strip.
 The cocycle is given by an integral of a form which includes derivatives
 over the period  (over the  elliptic nome in case of $\ap$) instead
 of  derivatives over spectral parameter.
 These algebras are isomorphic for
different strips.
 In contrast to \cite{RS} and
\cite{U}, there are natural isotropic subalgebras
  only
in the limit $\eta\to 0$ ($p\to 0$ in case of $\ap$).
Nevertheless at $c=0$
we have Lie bialgebra which includes according to Sokhotsky formulas
the usual currents on the line (on the circle in case of $\ap$).

On the other hand, we can describe the Lie algebra $\ael$
 in terms of Fourier harmonics of the generating functions for $\ael$.
 This language is a variant of the usual description of a current algebra
 in terms of Fourier modes
 for the case of vanishing at infinity currents on the real line.
The structural constants of the algebras $\ap$ or $\ael$ become standard
(do not depend on $p$ or $\eta$) in terms of Fourier components of
the generating functions, but the cobracket has a nontrivial form.

We would like to emphasize the ideology for the description of the algebra
 $\ael$ which we keep throughout the paper. The algebra  $\ael$
 cannot be described in terms of a discrete basis like Taylor
 coefficients of the analytical
functions. In order to define an analog of a basis
 we are forced to use the language of Fourier harmonics on an open line.
Since we get in this way continuously many generators we should specify which
integrals over them belong to the algebra. Equivalently, we fix a region
 for a spectral parameter where the generating functions are meromorphic.
 In representation theory it means that their matrix coefficients
 remain meromorphic in this region. Then we define an analytical
 continuation of the generating functions to larger domains of the spectral
 parameter. For the definition of this analytical continuation, as well as for
 the definition of the central extension, we need the language of generalized
function on a strip.
 Our construction of the algebra $\ael$ includes essentially all this stuff.

The ideology has been borrowed from \cite{KLP} and the main motivation
 of writing this paper was to clarify this ideology. We also present
an example
 of an infinite-dimensional representation of the algebra $\ael$ in the 
 Fock space.
 One can see that this representation requires  a description which differs
 from the usual construction of the
integrable representations of $\widehat{sl_2}$.

\setcounter{equation}{0}
\section{The Lie Algebra $\aelsect$ of the Generalized Automorphic\newline
 Functions on a Strip}

\subsection{The algebra $\aelcsseco$ and its central extension}

Let $e,f,h$ be standard generators of Lie algebra $\frak{sl}_2$:
$$
[h,e]=2e,\quad [h,f]=-2f,\quad [e,f]=h\ ,
$$
and $\eta=1/\xi\in\RR$ be a positive real parameter. Consider the meromorphic
$\frak{sl}_2$-valued functions of $z\in\CC$ with the period $2i/\eta$
which satisfy the following asymptotical and
automorphic conditions:
\bea
&(\iota)&  e(z)=-e(z+i/\eta),\quad
f(z)=-f(z+i/\eta),\quad h(z)=h(z+i/\eta),\nn\\
&(\iota\iota)& |e(x+
iy)|\stackreb{<}{x\to\pm\infty} C\ee^{-\pi\eta|x|},\quad |f(x+
iy)|\stackreb{<}{x\to\pm\infty} C\ee^{-\pi\eta|x|},\quad |h(x+
iy)|\stackreb{<}{x\to\pm\infty} C\ .  \nn \eea
One can verify that the
functions \bea e(z)&=&e\otimes {\cth^n\pi\eta(z-a)\over\sh\,
\pi\eta(z-a)},\quad
n\geq0, \quad a\in\CC\ ,\label{ez}\\ f(z)&=&f\otimes
{\cth^n\pi\eta(z-b)\over\sh\,\pi\eta(z-b)},\quad n\geq0,
\quad b\in\CC\ ,\label{fz}\\
h(z)&=&h\otimes \cth^n\pi\eta(z-c),\quad n\geq0,
\quad c\in\CC\ \label{hz}
\eea
satisfy the conditions $(\iota)$ and $(\iota\iota)$
and their finite linear combinations
form a Lie algebra denoted by $\aelt$.

The algebra $\aelt$ can be described by means of the generating
functions depending on the generating parameter $u$:
\beq
\tle_\eta(u)=e\otimes {i\pi\eta\over\sh\,\pi\eta(z-u)},\quad
\tlh_\eta(u)=h\otimes i\pi\eta\,\cth\,\pi\eta(z-u),\quad
\tlf_\eta(u)=f\otimes {i\pi\eta\over\sh\,\pi\eta(z-u)}\ ,
\label{gener}
\eeq
which satisfy the commutation relations
\bea
{[}\tlh_\eta(u_1), \tle_\eta(u_2){]}&=& 2i\pi\eta\,\cth\,\pi\eta(u_1-u_2)
\tle_\eta(u_2)-
{2i\pi\eta\over \sh\,\pi\eta (u_1-u_2)} \tle_\eta(u_1)\ ,\label{the}\\
{[}\tlh_\eta(u_1), \tlf_\eta(u_2){]}&=& -2i\pi\eta\,\cth\,\pi\eta(u_1-u_2)
\tlf_\eta(u_2)+
{2i\pi\eta\over \sh\,\pi\eta (u_1-u_2)} \tlf_\eta(u_1)\ ,\label{thf}\\
{[}\tle_\eta(u_1), \tlf_\eta(u_2){]}&=& {i\pi\eta\over \sh\,\pi\eta (u_1-u_2) }
\left(\tlh_\eta(u_1)-\tlh_\eta(u_2)\right)\ .\label{thh}
\eea

The Lie algebras $\aelt$ are isomorphic for different (finite) values of parameter $\eta$. An isomorphism
\beq
T_{\eta,\eta'}: \aelt\to
{\tilde{\frak{a}}_{\eta'}^0(\widehat{\frak{sl}_2})}
\label{iso}
\eeq
comes from the gauge transformation
\beq
a(u)\to \lambda a(\lambda u) \ .
\label{gauge}
\eeq
In terms of generating functions \r{gener} the isomorphism \r{iso} looks
 as follows:
\beq
\tle_{\eta'}(u)=\frac{\eta'}{\eta}\tle_\eta\left(\frac{\eta'}{\eta}u\right)
\ , \quad
\tlh_{\eta'}(u)=\frac{\eta'}{\eta}\tlh_\eta\left(\frac{\eta'}{\eta}u\right)
\ , \quad
\tlf_{\eta'}(u)=\frac{\eta'}{\eta}\tlf_\eta\left(\frac{\eta'}{\eta}u\right)
\ .
\label{gauge1}
\eeq

Let $\Pi^+\subset\CC$ be a strip $-1/\eta<\Im z<0$.
 In the sequel we need an interpretation of
the elements of the algebra
$\aelt$  as of generalized  $\frak{sl}_2$-valued function
on the strip $\Pi^+$. This description reads as follows.
Let $K$ be the space of basic functions $s(z)$ analytical
in  a strip $\Pi^+$,
 continuous in the closure $\overline{\Pi}^+$ of $\Pi^+$
and decreasing in the closed strip $\overline{\Pi}^+$
faster than some exponential function:
$$
|s(x+iy)|\stackreb{<}{x\to\pm\infty}C\ee^{-\alpha|x|}\ ,\qquad \alpha > 0\ .
$$

We treate $sl_2$-valued functions \r{ez}--\r{hz} as $sl_2$-valued functionals
 on the space $K$ and denote them by $e_+(z), h_+(z), f_+(z)$.
A pairing $(a_+,s)\in\frak{sl}_2$, $a\in\aelt$, $s\in K$
is defined as
\beq
\intt dx\  a(x)\, s(x)\ .
\label{pairing}
\eeq

Let us denote the described space of generalized functions as
$\aelo$. It inherits the structure of Lie algebra $\aelt$. As before, we
 compose the distributions into generating functions which we denote now by
$e_+(u), h_+(u)$, and $ f_+(u)$:
\beq
e_+(u)=e\otimes {i\pi\eta\over\sh\,\pi\eta(z-u)},\quad
h_+(u)=h\otimes i\pi\eta\,\cth\,\pi\eta(z-u),\quad
f_+(u)=f\otimes {i\pi\eta\over\sh\,\pi\eta(z-u)}\ .
\label{gener1}
\eeq
An index $+$ reminds now that the contour of integration in the pairing
 \r{pairing} lies above the pole $u$. For a fixed $u \in \Pi^+$ these
 generating functions are distributions from the space $\aelo$.

The algebra $\aelo$ is a Lie bialgebra. The cobracket $\delta$ is given by
the relations
\bea
\delta(e_+(u))&=&h_+(u)\wedge e_+(u)\ ,\nn\\
\delta(f_+(u))&=&f_+(u)\wedge h_+(u)\ ,\nn\\
\delta(h_+(u))&=&2e_+(u)\wedge f_+(u)\ .\label{6}
\eea
The algebra $\aelo$ possesses an invariant scalar product $\langle\, ,\,
\rangle$:
\beq
\langle e_+(u_1),f_+(u_2)\rangle = {\pi\eta^2 (u_1-u_2)
\over\sh\,\pi\eta (u_1-u_2)},\quad
\langle h_+(u_1),h_+(u_2)\rangle =2\pi\eta^2 (u_1-u_2)\,
\cth\,\pi\eta(u_1-u_2)\ .
\label{6a}
\eeq
Note that this scalar product  differs from the one used in \cite{RS}.
We will specify below the subalgebras of $\aelo$ which become
isotropic in the rational limit $\eta\to0$.

The Lie algebra $\aelo$ of generalized automorphic $\frak{sl}_2$-valued
functions on the strip  admits a central extension.
 It can be
 defined by a ``strange'' cocycle
\beq
B(x\otimes \varphi(z),y\otimes \psi(z))={\eta^2\over 4\pi}
\int_{\partial \Pi}dz \left({d\psi(z)\over d\eta}\varphi(z)-
\psi(z){d\varphi(z)\over d\eta}
\right) \ \langle x,y\rangle\ ,
\label{3.2}
\eeq
where $\langle\ ,\ \rangle$ is the Killing form, $x,y\in\frak{sl}_2$.
The integration in \r{3.2} goes anticlockwise along the boundary
$\partial \Pi$ of the strip $\Pi$
which consists of the real axis
and of the line $\Im z=-1/\eta$.
The commutation relations of the extended algebra $\ael$ are
($u=u_1-u_2$)
\bea
{[}h_+(u_1), e_+(u_2){]}&=& 2i\pi\eta\,\cth(\pi\eta u) e_+(u_2)-
{2i\pi\eta\over \sh\,\pi\eta u} e_+(u_1)\ ,\nn\\
{[}h_+(u_1), f_+(u_2){]}&=& -2i\pi\eta\,\cth(\pi\eta u) f_+(u_2)+
{2i\pi\eta\over \sh\,\pi\eta u} f_+(u_1)\ ,\nn\\
{[}e_+(u_1), f_+(u_2){]}&=& {i\pi\eta\over \sh\,\pi\eta u }
\left(h_+(u_1)-h_+(u_2)\right)+ ic\pi\eta^2\left(
{\pi\eta u\,\ch\,\pi\eta u\over\sh^2\,\pi\eta u}-
{1\over\sh\,\pi\eta u}\right),\nn\\
{[}h_+(u_1), h_+(u_2){]}&=& 2ic\pi\eta^2\left(
{\pi\eta u\over\sh^2\,\pi\eta u}-\cth\,\pi\eta u\right).\label{2.3}
\eea
Unfortunately, we do not know how to extend the Lie coalgebra structure
 \r{6} from the algebra $\aelo$ to its central extension
 $\ael$. We return to this point later in section 4.

The isomorphism \r{iso} induced by the gauge transformation \r{gauge}
 has a natural prolongation to the central extended algebra $\ael$. 
 It preserves
 the pairing \r{6a} and multiply the cobracket $\delta: \aelo \to \aelo \ot
\aelo$ by a factor $\eta/\eta'$. So we have actually the unique algebraic
 structure realized in different spaces of distributions.

\subsection{The analytical continuation and the Sokhotsky's formulas}

{} For a periodic function (over variable $z$) of the type
$$a(u)=
e\otimes {i\pi\eta\over\sh\,\pi\eta(z-u)}\quad {\rm or}
\quad a(u)=h\otimes i\pi\eta\,\cth\,\pi\eta(z-u)\quad
{\rm or}\quad a(u)=f\otimes {i\pi\eta\over\sh\,\pi\eta(z-u)}$$
 let us denote by $a_-(u)$ the distribution whose value on a basic
 function $s(z) \in K$ is given by the integral
$$\int_{\Gamma} dz\ a(z)s(z)$$
where the contour $\Gamma$ is a line parallel to the real line such that
 the pole $u$ is the first pole of $a(z)$ located above the contour
 $\Gamma$. Then for $u \in \Pi^- = \{0<\Im z< 1/\eta\}$ we have by definition
\beq
e_-(u)=-e_+(u-i/\eta),\quad
h_-(u)=h_+(u-i/\eta),\quad
f_-(u)=-f_+(u-i/\eta)\ .
\label{conn}
\eeq
The distributions (over the variable $z$) $e_\pm(u)$, $f_\pm(u)$ and
$h_\pm(u)$ admit analytical continuations over the parameter $u$. For instance,
 the analytical continuation of $e_-(u), f_-(u)$ and of $h_-(u)$
 to the region of parameter $u \in \Pi^+$ are the distributions of the type
$a_-(z)$ defined as
$$
(a_-,s)=
\int_{-i/\eta-\infty}^{-i/\eta+\infty} dz\
a(z)\, s(z)\ .
$$
 Due to the relations \r{conn} this analytical continuation preserves the
 structure of Lie algebra.
 The analytical continuation of the commutation relations
\r{2.3} yields:
\bea
{[}e_\pm(u_1), f_\pm(u_2){]}&=& {i\pi\eta\over \sh\,\pi\eta (u_1-u_2) }
\left(h_\pm(u_1)-h_\pm(u_2)\right)+B(e_\pm(u_1),f_\pm(u_2))\cdot c\ ,\nn\\
{[}h_\pm(u_1),h_\pm(u_2){]}&=& B(h_\pm(u_1),h_\pm(u_2))\cdot c\ ,\nn\\
{[}h_\pm(u_1), e_\pm(u_2){]}&=& 2i\pi\eta\,\cth\,\pi\eta(u_1-u_2)e_\pm(u_2)-
{2i\pi\eta\over \sh\,\pi\eta (u_1-u_2)} e_\pm(u_1)\ ,\nn\\
{[}h_\pm(u_1), f_\pm(u_2){]}&=&-2i\pi\eta\,\cth\,\pi\eta(u_1-u_2)f_\pm(u_2)+
{2i\pi\eta\over \sh\,\pi\eta (u_1-u_2)} f_\pm(u_1)\ ,\nn\\
{[}e_\pm(u_1), f_\mp(u_2){]}&=& {i\pi\eta\over \sh\,\pi\eta (u_1-u_2) }
\left(h_\mp(u_1)-h_\pm(u_2)\right)+B(e_\pm(u_1),f_\mp(u_2))\cdot c\ ,\nn\\
{[}h_\pm(u_1),h_\mp(u_2){]}&=& B(h_\pm(u_1),h_\pm(u_2))\cdot c\ ,\nn\\
{[}h_\pm(u_1), e_\mp(u_2){]}&=& 2i\pi\eta\,\cth\,\pi\eta(u_1-u_2)e_\mp(u_2)-
{2i\pi\eta\over \sh\,\pi\eta (u_1-u_2)} e_\pm(u_1)\ ,\nn\\
{[}h_\pm(u_1), f_\mp(u_2){]}&=&-2i\pi\eta\,\cth\,\pi\eta(u_1-u_2)f_\mp(u_2)+
{2i\pi\eta\over \sh\,\pi\eta (u_1-u_2)} f_\pm(u_1)\ ,\label{3.1}
\eea
where  $(u=u_1-u_2)$
\bea
B(e_\pm(u_1),f_\pm(u_2))
&=& i\pi\eta^2\left(
{\pi\eta u\,\ch\,\pi\eta u\over\sh^2\,\pi\eta u}
-{1\over\sh\,\pi\eta u}\right)
\ ,\nn\\
B(e_\pm(u_1),f_\mp(u_2))
&=& i\pi\eta^2\left(
{\pi\eta(u\pm i/\eta)\,\ch\,\pi\eta u\over\sh^2\,\pi\eta u}
-{1\over\sh\,\pi\eta u}\right)
\ ,\nn\\
B(h_\pm(u_1),h_\pm(u_2))
&=& 2i\pi\eta^2\left(
{\pi\eta u\over\sh^2\,\pi\eta u}
-\cth\,\pi\eta u\right)
\ ,\nn\\
B(h_\pm(u_1),h_\mp(u_2))
&=& 2i\pi\eta^2\left(
{\pi\eta(u\pm i/\eta)\over\sh^2\,\pi\eta u}
-\cth\,\pi\eta u\right)
\ .\nn
\eea
The pairing \r{6a} also admits an analytical continuation:
\bea
\langle e_\pm(u_1),f_\mp(u_2)\rangle &=& {\pi\eta^2 (u_1-u_2\pm i/\eta)
\over\sh\,\pi\eta (u_1-u_2)}\ ,\nn\\
\langle h_\pm(u_1),h_\mp(u_2)\rangle &=&2\pi\eta^2 (u_1-u_2\pm i/\eta)\,
\cth\,\pi\eta(u_1-u_2)\ .\label{6b}
\eea

The definition of the distributions
$e_\pm(u)$, $f_\pm(u)$, $h_\pm(u)$
 gives possibility to apply Sokhotsky's formulas. They can be written
 as follows:
\beq e(u)=e_+(u)-e_-(u)\ , \quad f(u)=f_+(u)-f_-(u),\quad
h(u)=h_+(u)-h_-(u)\ ,
\label{4a}
\eeq
where
\beq
e(u)=2\pi e\ \otimes\ \delta(u)\   \ ,\qquad
f(u)=2\pi f\ \otimes\ \delta(u)\ ,\qquad
h(u)=2\pi h\ \otimes\ \delta(u)\ .
\label{4c}
\eeq
On the other hand, we also have  the relations
\bea
e_\pm(u)&=&{i\eta\over2}\int_{\Gamma_\pm}dz\ {e(z)\over\sh\,\pi\eta(z-u)}\ ,
\nn\\
f_\pm(u)&=&{i\eta\over2}\int_{\Gamma_\pm}dz\ {f(z)\over\sh\,\pi\eta(z-u)}\ ,
\nn\\
h_\pm(u)&=&{i\eta\over2}\int_{\Gamma_\pm}dz\ {h(z)\ \cth\,\pi\eta(z-u)}\ ,
\label{7}
\eea
where the countor $\Gamma_+$ is a line parallel to the real axis and lying
above the point $u$ and $\Gamma_-$ is also a line parallel to the real
axis but the point $u$ is above  it.

The relations \r{4a}--\r{7} show that the algebra of formal $sl_2$
valued currents on the line is embedded into the analytical continuation of
 the Lie algebra $\ael$:

\bea
{[}h(u),e(v){]}&=&2\delta(u-v)e(v)\ ,\nn\\
{[}h(u),f(v){]}&=&-2\delta(u-v)f(v)\ ,\nn\\
{[}e(u),f(v){]}&=&\delta(u-v)h(v)+c\cdot\delta'(u-v)\ ,\nn\\
{[}h(u),h(v){]}&=&2c\cdot\delta'(u-v)\ .
\label{4b}
\eea

\setcounter{equation}{0}
\section{The Lie Algebra $\aelsect$ in  Terms of  Fourier Harmonics}

Let $\he_\la, \hf_\la$ and $\hf_\la$, $\la\in\RR$
  be the symbols satisfying the following relations
\bea
{[}\hh_\la,\he_\mu{]}&=&2\he_{\la+\mu}\ ,\nn\\
{[}\hh_\la,\hf_\mu{]}&=&-2\hf_{\la+\mu}\ ,\nn\\
{[}\he_\la,\hf_\mu{]}&=&\hh_{\la+\mu}+c\cdot\la\delta(\la+\mu)\ ,\nn\\
{[}\hh_\la,\hh_\mu{]}&=&2c\cdot\la\delta(\la+\mu)\ .
\label{9}
\eea

Let $\overline{\frak{a}}_\eta$ be a vector space which consists of the
 expressions of the type
\beq
\int_{-\infty}^\infty d\la\
{\he_\la g(\la) \over 1+\ee^{\la/\eta}},\quad
\int_{-\infty}^\infty d\la\
{\hf_\la g'(\la) \over 1+\ee^{\la/\eta}},\quad
{}\vpint d\la\
{\hh_\la g''(\la) \over 1-\ee^{\la/\eta}}\ ,
\label{11}
\eeq
where $g(\la)$, $g'(\la)$, $g''(\la)$ are quasi-polynomials
\beq
g(\la)=\sum_jP_j(\la)\ee^{i\la u_j},\quad
g'(\la)=\sum_jQ_j(\la)\ee^{i\la u_j},\quad
g''(\la)=\sum_jR_j(\la)\ee^{i\la u_j}\ ,\label{12}
\eeq
$u_j\in\Pi^+$ and $P_j(\la)$, $Q_j(\la)$, $R_j(\la)$ are polynomials.

We state that the brackets \r{9} define 
a Lie algebra structure on the space
$\overline{\frak{a}}_\eta$ which is isomorphic to
$\ael$.\\
The isomorphism can be established in the following way. We realize the
symbols
 $\he_\la$, $\hf_\la$
and $\hh_\la$ as the following functionals on the space $K$
 of basic functions $s(z)$:
\bea
(\he_\la,f(z))&=&e\ot\intt dz\,s(z)e^{-i\la z}\ ,
\label{el}\\
(\hf_\la,f(z))&=&f\ot\intt dz\,s(z)e^{-i\la z}\ ,
\label{hl}\\
(\hh_\la,f(z))&=&h\ot\intt dz\,s(z)e^{-i\la z}\ .
\label{fl}
\eea

Then the expressions \r{11} are in a one-to-one correspondence 
with distributions
from $\ael$ defined by the functions \r{ez}--\r{hz}. Since both brackets
 \r{the}--\r{thh} and \r{9} for $c= 0$ are pointwise brackets of $sl_2$ valued
 functions, they coincide for $c=0$. One can easily check that the cocycles
in \r{2.3} and in \r{9} also coincide.

In the language of Fourier harmonics we can define a natural extension of the
 algebra $\ael$. This extension consists of the following expressions
\beq
\int_{-\infty}^\infty d\la\
\he_\la  \tilde g(\la),\quad
\int_{-\infty}^\infty d\la\
\hf_\la  \tilde g'(\la),\quad
{}\vpint d\la\
\hh_\la  \tilde g''(\la) \ ,
\label{ggg}
\eeq
where these functions satisfy the following conditions:
\bea
\tilde g(\la) \;{\rm and}\;  \tilde g'(\la)
&\ \ \mbox{are analytical in the strip}\
&-\pi\eta<\Im(\la)<\pi\eta\  ,\nn\\
\tilde g''(\mu)\qquad
&\ \ \mbox{is analytical in the strip}\
&-2\pi\eta<\Im\mu<2\pi\eta\ ,\nn
\eea
except the point $\mu=0$ where the function $\tilde g''(\mu)$ has a simple
pole. Besides, the functions $\tilde g(\la)$,
$\tilde g'(\la)$, $\tilde g''(\la)$ should decrease
 faster than some exponent  when
$\Re\la\to\pm\infty$:
$$\tilde g(\la)< C\ee^{-\alpha|\Re\la|}\ ,\qquad
 \tilde g'(\la)< C\ee^{-\beta|\Re\la|}\ , \qquad
 \tilde g''(\la)< C\ee^{-\gamma|\Re\la|}\ $$
for some $\alpha, \beta ,\gamma >0.$

This extended algebra is well-defined and the arguments are
the same as in
\cite{KLP}. We do not reserve  any special notation for the extended algebra.
 It plays the same role as $\ael$.

All the structures of the Lie algebra $\ael$ can be reformulated
 in the language of Fourier harmonics. The invariant scalar product \r{6a} now
 looks standard:
\beq
\langle\he_\la,\hf_\mu\rangle=\delta(\la+\mu),\quad
\langle\hh_\la,\hh_\mu\rangle=2\delta(\la+\mu)\
\label{can-pair}
\eeq
and the cobracket for level $0$ algebra $\aelo$ is:
\bea
\delta \he_\la&=&-{1\over2}\vpint d\tau\ \hh_\tau\wedge\he_{\la-\tau}
\left(\cth\,\tau/2\eta+\th\,(\la-\tau)/2\eta \right)
\ ,\nn\\
\delta \hf_\la&=&-{1\over2}\vpint d\tau\ \hf_{\la-\tau}\wedge\hh_\tau
\left(\cth\,\tau/2\eta+\th\,(\la-\tau)/2\eta \right)
\ ,\nn\\
\delta \hh_\la&=&-\vpint d\tau\ \he_\tau\wedge\hf_{\la-\tau}
\left(\th\,\tau/2\eta+\th\,(\la-\tau)/2\eta \right)
 \ .\label{3.3}
\eea
The generating functions $e_\pm(u), f_\pm(u)$ and $h_\pm(u)$ are treated
 now as generating integrals for $\he_\la, \hf_\la, \hh_\la$.
As it follows from
 \r{el}--\r{fl} and from \r{gener1},
\bea
e_\pm(u)&=&\pm \int_{-\infty}^\infty d\la\ \ee^{i\la u}\
{\he_\la \over 1+\ee^{\pm\la/\eta}}\ ,\nn\\
f_\pm(u)&=&\pm \int_{-\infty}^\infty d\la\ \ee^{i\la u}\
{\hf_\la \over 1+\ee^{\pm\la/\eta}}\ ,\nn\\
h_\pm(u)&=&\pm\  \vpint d\la\ \ee^{i\la u}\
{\hh_\la \over 1-\ee^{\pm\la/\eta}}\ .\label{8a}
\eea
and the currents  $e(u)$, $f(u)$ and $h(u)$ are total Fourier images of
 $\he_\la$, $\hf_\la$, $\hh_\la$:
\beq
e(u)=\intt d\la\ \he_\la \ee^{i\la u},\quad
f(u)=\intt d\la\ \hf_\la \ee^{i\la u},\quad
h(u)=\intt d\la\ \hh_\la \ee^{i\la u}\ .
\label{8}
\eeq
\bigskip

\setcounter{equation}{0}
\section{The Representation Theory}

As we already mentioned in Introduction the matrix elements
of the generating functions $h_+(u)$,  $e_+(u)$,  $f_+(u)$ become
the meromorphic functions in the representations of the algebra
$\ael$. Let us define representations of this algebra in the
description given by \r{ggg}. Assign to each element defined by
\r{ggg}, for example, $\intt d\la g(\la)\he_\la$ the operator valued
function:
$$
\tilde e(u)=\intt d\la g(\la) \he_\la \ee^{i\la u}\ .
$$
We will say that a representation of the algebra $\ael$ is well-defined
if the operator valued function $\tilde e(u)$ becomes a
meromorphic function in the variable $u\in\CC$ in some neighbourhood of zero.

\subsection{Representation of the algebra $\aelcssect$ at the level 1}

The goal of this subsection is to construct an infinite-dimentional
representation of the algebra $\aelc$ at level $c=1$.
 For a description of this representation we need
 a definition of the Fock space generated by continuous family of free
bosons. We borrow this definition from \cite{KLP}.

Let $a_\la$, $\la\in\RR$, $\la\neq 0$ be free bosons which satisfy the
 commutation relations
$$
[a_\la, a_\mu]=a(\la)\,\delta(\la+\mu)\ ,\qquad a(\la)={\la\over2}\ .
$$

 We define a  (right) Fock space $\H_{a(\la)}$ as follows.
$\H_{a(\la)}$ is generated as a vector
space
by the expressions
$$
\mint f_1(\la_1) a_{\la_1} d\la_1\ldots \mint f_n(\la_n) a_{\la_n} d\la_n
\ \rvac\  ,
$$
 where the functions $f_i(\la)$ satisfy the condition
$$
f_i(\la)< C\, \ee^{\epsilon\la}, \qquad \la \to -\infty\ ,
$$
for $\epsilon>0$ and
$f_i(\la)$ are analytical functions in a neighbourhood of
 ${\RR_+}$ except $\la=0$, where they have a simple pole.

The left Fock space $\H^*_{a(\la)}$ is generated by the expressions
$$
\lvac \nint g_1(\la_1) a_{\la_1} d\la_1\ldots \nint g_n(\la_n) a_{\la_n}
d\la_n\ ,
$$
where the functions $g_i(\la)$ satisfy the conditions
$$
g_i(\la)< C\, \ee^{-\epsilon\la}, \qquad \la \to +\infty\ ,
$$
and $g_i(t)$ are analytical functions in a neighbourhood of
 ${\RR_-}$ except $\la=0$, where they also have a simple pole.

The pairing $(,):$ $\H_{a(\la )}^*\otimes \H_{a(\la )} \to \CC$ is
uniquely defined by
the
 following prescriptions:
\bea
 &(i)&\quad  (\langle\mbox{vac}|, \rvac) =1\ ,\nn\\
 &(ii)& \quad (\lvac \nint d\la \ g(\la)a_\la\ ,
\mint d\mu\  f(\mu)a_\mu\   \rvac) =
\int_{\tilde C} {d\la\,\ln(-\la)\over 2\pi i}
  g(\la)f(-\la)a(\la)\ ,\nn\\
&(iii)& \quad \mbox{the Wick theorem}\ \nn
\eea
and contour $\tilde C$ shown on the Fig.~1.
\bigskip

\unitlength 1.00mm
\linethickness{0.4pt}
\begin{picture}(121.00,20.00)
\put(17.00,15.00){\makebox(0,0)[cc]{0}}
\put(20.00,15.00){\makebox(0,0)[cc]{$\bullet$}}
\put(132.00,15.00){\makebox(0,0)[cc]{$+\infty$}}
\put(20.00,15.00){\line(1,0){100.33}}
\put(40.00,10.00){\line(1,0){80.33}}
\put(120.00,20.00){\line(-1,0){100.00}}
\put(30.00,5.00){\makebox(0,0)[cc]{Fig.~1.}}
\put(121.00,10.00){\vector(1,0){0.2}}
\put(100.00,10.00){\line(1,0){21.00}}
\put(20.00,10.00){\line(1,0){22.00}}
\put(20.67,10.00){\line(1,0){22.00}}
\put(20.17,15.00){\oval(15.00,10.00)[l]}
\end{picture}
\smallskip

Let the vacuums $\lvac$ and $\rvac$ satisfy
 the conditions
$$
a_\la\rvac =0,\quad \la>0,\qquad \lvac a_\la =0, \quad \la<0\ ,
$$
 and $f(\la)$ be  a function analytical in some neighbourhood of the real line
 with possible simple pole at $\la=0$ and which has the
following asymptotical behaviour:
$$
\qquad f(\la)< e^{-\epsilon|t|}
, \qquad \la \to \pm\infty
$$
for some $\epsilon>0$.
Then, by definition, an operator
$$F= {:}\exp \sk{\int_{-\infty}^{+\infty}d\la\ f(\la)a_\la} {:}$$
acts on the right Fock space $\H_{a(\la )}$ as follows.
$F=F_-F_+$, where
$$F_-=\exp \sk{\int_{-\infty}^{0}d\la\ f(\la)a_\la} \quad
\mbox{and}\quad
F_+=
\lim{\epsilon\to 0} \epsilon^{
\epsilon f(\epsilon)a_{\epsilon}}
\exp\sk{\int_{\epsilon}^{\infty}d\la\ f(\la)a_\la}\ .
$$
An action of operator $F$ on the left Fock space $\H_{a(\la )}^*$ is
 defined via another decomposition: $F=\tilde{F}_-\tilde{F}_+$, where
$$\tilde{F}_+=\exp \sk{\nint d\la\  f(\la)a_\la}\quad
\mbox{and}\quad
\tilde{F}_-=
\lim{\epsilon\to 0} \epsilon^{
\epsilon f(-\epsilon)a_{-\epsilon}}
\exp\sk{\int^{-\epsilon}_{-\infty}d\la\ f(\la)a_\la}.
$$
 These definitions imply
that the defined above actions of the operator
$$
F= :\exp \sk{\int_{-\infty}^{+\infty}d\la\ f(\la)a_\la}{:}
$$
on the Fock spaces $\H_{a(\la)}$ and $\H_{a(\la)}^*$ are adjoint and the
product of normally ordered operators satisfy the property \cite{MJ2}
\bea
&{:}\exp\left(\intt d\la\ g_1(\la)\,a_\la\right){:}\ \cdot\
{:}\exp\left(\intt d\mu\ g_2(\mu)\,a_\mu\right){:} = \nn\\
&\quad =\exp\left(\int_{\tilde C}{d\la\,\ln(-\la)\over2\pi i}\
a(\la)g_1(\la)g_2(-\la) \right)
{:}\exp\left(\intt d\la\ (g_1(\la)+g_2(\la))\,a_\la\right){:}\ .
\label{normal}
\eea
\medskip

\noindent
{\bf Proposition.} {\it The generating functions} \r{last}
{\it satisfy the commutation relations} \r{3.1} {\it and
define a highest  weight, level $1$ representation of the algebra $\aelc$.}
\bea
e_\pm(u)&=&{i\eta\ee^\gamma\over2}
\int_{\Gamma_\pm} {dz\over\sh\,\pi\eta(z-u)}
{:}\exp\left(\intt d\la\ \ee^{i\la z}{2a_\la\over \la}\right){:}\ ,
\nn\\
f_\pm(u)&=&{i\eta\ee^\gamma\over2}
\int_{\Gamma_\pm} {dz\over\sh\,\pi\eta(z-u)}
{:}\exp\left(-\intt d\la\ \ee^{i\la z}{2a_\la\over \la}\right){:}\ ,
\label{last}\\
h_\pm(u)&=&\pm\  \vpint d\la\ \ee^{i\la u}\
{2a_\la \over 1-\ee^{\pm\la/\eta}}\ . \nn
\eea

The contours $\Gamma_\pm$ in \r{last} are the same as in \r{7} and
$\gamma$ is Euler constant. The action of the total currents $e(u)$, $f(u)$ and
$h(u)$ has a simple form
$$e(u)=\ee^\gamma
{:}\exp\left(\intt d\la\ \ee^{i\la u}{2a_\la\over \la}\right){:}\ ,
\qquad
f(u)=\ee^\gamma
{:}\exp\left(-\intt d\la\ \ee^{i\la u }{2a_\la\over \la}\right){:}\ ,$$
 $$
h(u)= 2\intt d\la\ a_\la\ee^{i\la u}\ .$$

\setcounter{equation}{0}
\section{The Algebra $\Aelsect$ and its Classical Limit}

\subsection{The definition of $\Aelssect$}

The algebra $\Ael$ is a scaling limit of the elliptic algebra
$\Apq$ \cite{MJ} at $p,q\to1$. This algebra has been investigated
in \cite{KLP}. $\Ael$ is generated by the central element
$c$ and by certain integrals over  Fourier harmonics $\he_\la, \hf_\la,
 \hh_\la$ for $c =0$ and  $\he_\la, \hf_\la,
 \hat{t}_\la$ for $c \not =0$ (see \cite{KLP})
of the matrix elements of the $L$-operator $L^+(u,\eta)$ or its analytical
continuation $L^-(u,\eta)$:
\beq
L^-(u,\eta)=
\sigma_z   L^+\left(u-i(1/\eta+\tih c/2),\eta\right)  \sigma_z\ .
\label{2.0}
\eeq
Operator $L^+(u,\eta)$ satisfies the relations
\begin{eqnarray}
R^+(u_1-u_2,\eta')L^+_1(u_1,\eta)L^+_2(u_2,\eta)&=&
L^+_2(u_2,\eta)L^+_1(u_1,\eta) R^+(u_1-u_2,\eta) \  , \label{2.1}\\
\qdet L(u,\eta)&=&1\ ,\nn
\eea
where
$$\eta'={\eta\over 1+\eta c\tih},\quad c\tih >0\ .$$
The last inequality means that in the
  representations   the
central element
$c$ is equal to some number $c$ such that
$\tih c>0$.
The $R$-matrix in \r{2.1} reads as
\begin{eqnarray}
R^+(u,\eta)&=&\tau^+(u)  R (u,\eta),\quad
R(u,\eta)\ =\ \varrho(u,\eta)\overline R (u,\eta)\ ,
\nn\\
\overline R (u,\eta)&=&
\left(
\begin{array}{cccc}
1&0&0&0\\  0&b (u,\eta)&c (u,\eta)&0\\
0&c (u,\eta)&b (u,\eta)&0\\  0&0&0&1
\end{array}
\right)      \ ,     \nn\\
\varrho (u,\eta)&=&{\Ga{\tih\eta}\Ga{1+i\eta u}\over
        \Ga{\tih\eta+i\eta u}}
\prod_{p=1}^\infty
{R _p(u,\eta) R _p(i\tih- u,\eta)
\over
R _p(0,\eta) R _p(i\tih,\eta)}\ , \nn\\
R _p(u,\eta)&=&{\Ga{2p\tih\eta+i\eta u}
               \Ga{1+2p\tih\eta+i\eta u}\over
        \Ga{(2p+1)\tih\eta+i\eta u}
        \Ga{1+(2p-1)\tih\eta+i\eta u}}\ ,\nn\\
b(u,\eta)\ &=&\
{\sh\,\pi\eta u\over\sh\,\pi\eta(u-i\tih)}\ ,\quad
c(u,\eta)\ =\
{-\sh\,i\pi\eta\tih \over\sh\,\pi\eta(u-i\tih)}\ ,\quad
\tau^+(u)=\cth\sk{{\pi u\over2\tih}}    \label{R-mat}
\end{eqnarray}
and the quantum determinant is
\bea
\qdet L(z,\eta)
= L_{11}(z-i\tih,\eta)L_{22}(z,\eta)-
L_{12}(z-i\tih,\eta)L_{21}(z,\eta)\ .
\label{qdet}
\eea
Relations 
between 
$L^+(u)$ and
$L^-(u)$ can be obtained by means of the analytical continuations
of the relations \r{2.1}.

In terms of the coordinates  of the Gauss decomposition
\beq
L^+(u)
=\left(\begin{array}{cc} 1& f_+(u)\\0&1\end{array}\right)
\left(\begin{array}{cc}(k_+(u+i\tih))^{-1}&0\\0&k_+(u)\end{array}\right)
\left(\begin{array}{cc} 1&0\\ e_+(u)&1\end{array}\right)\ ,
\label{GL-univ}
\eeq
the relations \r{2.1} read as follows ($u=u_1-u_2$):

\bea
e_+(u_1)f_+(u_2)-f_+(u_2)e_+(u_1)&=&
{\sh\,i\pi\eta'\tih\over\sh\,\pi\eta'u}h_+(u_1)-
{\sh\,i\pi\eta\tih\over\sh\,\pi\eta u}\tilde h_+(u_2),\nn\\
\sh\,\pi\eta(u+i\tih)h_+(u_1)e_+(u_2)-
\sh\,\pi\eta(u-i\tih)e_+(u_2)h_+(u_1)&=&\sh(i\pi\eta\tih)
\{h_+(u_1),e_+(u_1)\},\nn\\
\sh\,\pi\eta'(u-i\tih)h_+(u_1)f_+(u_2)-
\sh\,\pi\eta'(u+i\tih)f_+(u_2)h_+(u_1)&=&-\sh(i\pi\eta'\tih)
\{h_+(u_1),f_+(u_1)\},\nn\\
\sh\,\pi\eta(u+i\tih)e_+(u_1)e_+(u_2)-
\sh\,\pi\eta(u-i\tih)e_+(u_2)e_+(u_1)&=&\sh(i\pi\eta\tih)
\left(e_+(u_1)^2+e_+(u_2)^2\right),\nn\\
\sh\,\pi\eta'(u-i\tih)f_+(u_1)f_+(u_2)-
\sh\,\pi\eta'(u+i\tih)f_+(u_2)f_+(u_1)&=&-\sh(i\pi\eta'\tih)
\left(f_+(u_1)^2+f_+(u_2)^2\right),\nn
\eea
\beq
\sh\,\pi\eta(u+i\tih)\sh\,\pi\eta'(u-i\tih)
h_+(u_1)h_+(u_2)
=h_+(u_2)h_+(u_1)
\sh\,\pi\eta'(u+i\tih)\sh\,\pi\eta(u-i\tih),\label{2.2}
\eeq
where
\bea
h_+(u)&=&k_+\left(u\right)^{-1}
k_+\left(u+i\tih\right)^{-1}\ , \nn\\
\tilde h_+(u)&=&k_+\left(u+i\tih\right)^{-1}
k_+\left(u\right)^{-1}= {\eta\over\eta'}
{\sin\,\pi\eta'\tih\over\sin\,\pi\eta\tih} h_+(u)
\ . \nn
\eea
The commutation relations \r{2.2} should be treated as the relations for the
 generating integrals $e_+(u)$, $f_+(u)$, $h_+(u)$ of elements  $\he_\la$,
$\hf_\la$,  $\hh_\la$.

The commutation relations for the operator $L^-(u)$ entries
can be obtained from \r{2.2} by the analytical continuation as it
was done in case of the commutation relations \r{3.1}.

The Hopf structure on $\Aelx$ is defined in a following sense.
 Consider the family of the algebras $\Aelx$ with fixed $\tih>0$  and
 variable $\xi=1/\eta$. Then the operations
$$
\Delta\,c=c_1+c_2=c\otimes1+1\otimes c\ ,$$
\beq
\Delta L^+_{ij}(u,\xi)=
\sum_{k=1}^2
L^+_{kj}(u+ic_2\tih/4,\xi) \otimes L^+_{ik}(u-ic_1\tih/4,\xi+\tih c_1)
\label{com-L*cmp}
\eeq
define a map
\beq
\Aelx\rightarrow \Aelx \otimes
{\cal A}_{\tih,1/(\xi+\tih c_1)}(\widehat{\frak{sl}_2})\label{map}
\eeq
which is coassociative and is
compatible with the commutation relations \r{2.1}.
The coproduct in terms of the generating functions $h_+(u)$,  $e_+(u)$,
$f_+(u)$ is given in \cite{KLP}.

\subsection{Classical limit $(\tih\to0)$ of the algebra $\Aelssect$}

Let $L^{\pm}(u,\eta) =1+\tih{\cal L}^\pm(u,\eta) +o(\tih^2)$. It is easy
to calculate that
\bea
\overline R(u,\eta)&=&
1+\tih r_0(u,\eta)+\tih (i\pi\eta\,\cth\,\pi\eta u)\ \id\ot\id +
\tih^2 r_1(u,\eta)  +o(\tih^2)\ ,\nn\\
\overline R(u,\eta')&=&
1+\tih r_0(u,\eta)+\tih (i\pi\eta\,\cth\,\pi\eta u)\ \id\ot\id +
\tih^2 r'_1(u,\eta)  +o(\tih^2)\ ,\nn\\
{\varrho(u,\eta')\over\varrho(u,\eta)}
&=& 1+\tih^2\varrho_0(u,\eta)+o(\tih^2)
\ ,\nn
\eea
where
\beq
r_0(u,\eta)=-i\pi\eta\left(
\begin{array}{cccc}
\cth\,\pi\eta u&0&0&\\ 0&0&(\sh\,\pi\eta u)^{-1}&0\\
0&(\sh\,\pi\eta u)^{-1}&0&0\\ 0&0&0&\cth\,\pi\eta u
\end{array}\right)
\label{r-clas}
\eeq
is a trigonometric solution to the classical Yang-Baxter
equations,
\bea
r_1(u,\eta)-r'_1(u,\eta)
&=&{i\pi\eta^2}\left(
\begin{array}{cccc}
0&0&0&0\\
0&\cth\,\pi\eta u-{\dis \pi\eta u\over\dis \sh^2\,\pi\eta u}
&{\dis \pi\eta u\,\ch\,\pi\eta u\over\dis \sh^2\,\pi\eta u}-{\dis1\over\dis
\sh\,\pi\eta u}&0\\
0&{\dis \pi\eta u\,\ch\,\pi\eta u\over\dis\sh^2\,\pi\eta u}-{\dis1\over\dis
\sh\,\pi\eta u}&
\cth\,\pi\eta u-{\dis \pi\eta u\over\dis \sh^2\,\pi\eta u}
&0\\
0&0&0&0
\end{array}\right) .\nn
\eea
and
$$
\varrho_0(u,\eta)=
{i\pi\eta^2\over2}\left(
\cth\,\pi\eta u-{\dis \pi\eta u\over\dis \sh^2\,\pi\eta u}\right).
$$
This expression can be found from the integral representation of the
factor $\varrho(u,\eta)$.

We obtain now from \r{2.1} the relation
\bea
{[}{\cal L}^+_1(u_1),{\cal L}^+_2(u_2){]}&=&{[}{\cal L}^+_1(u_1)+
{\cal L}^+_2(u_2),
r_0(u_1-u_2,\eta){]} +\nn\\
&+& \left(r_1(u_1-u_2,\eta) -r'_1(u_1-u_2,\eta)\right)c
-\varrho_0(u_1-u_2,\eta) c \cdot\id\otimes\id\ .
\label{LL}
\eea
These commutation relations can be found without calculation of the
expansion of the ratio of scalar factors
$\varrho(u,\eta')/\varrho(u,\eta)$. 
The role of the factor $\varrho_0(u,\eta)$ is to transform the matrix
in front of the central element to a 
traceless matrix. This condition fixes
the factor  $\varrho_0(u,\eta)$ uniquely.
Using the freedom
$$
r_0(u,\eta)\to \tilde r_0(u,\eta)=r_0(u,\eta)+\kappa(u)\cdot\id\ot\id
$$
let us make the new $r$-matrix $\tilde r_0(u,\eta)$ traceless.
Then the commutation relations \r{LL} can be written in the form
$$
{[}{\cal L}^+_1(u_1),{\cal L}^+_2(u_2){]}={[}{\cal L}^+_1(u_1)+
{\cal L}^+_2(u_2),
\tilde r_0(u_1-u_2,\eta){]} +
\eta^2{d\tilde r_0(u_1-u_2,\eta)\over d\eta}.
$$
We will use this observation in last section calculating the
classical limit of the quantum elliptic algebra $\Apq$.

Let $e_\pm(u)$, $f_\pm(u)$ and $h_\pm(u)$ be formal
 integrals of the symbols $\he_\la$, $\hf_\la$, $\hh_\la$
given by the formulas \r{8a} with the spectral parameter being
a complex number $u\in\CC$.
By the direct verification we can check that if the complex number
$u$ is inside the strip $\Pi^+$
then elements $e_+(u)$, $f_+(u)$, $h_+(u)$ belong to $\ael$. If the
complex number $u$ is inside the strip $\Pi^-$
then the elements
$e_-(u)$, $f_-(u)$, $h_-(u)$ also belong to $\ael$.
Thus we can treat the integrals
 $e_\pm(u)$, $f_\pm(u)$, $h_\pm(u)$
as generating functions of the elements of the algebra $\ael$, analytical
in the strips $\Pi^{\pm}$. We  can state the following
\smallskip

\noindent
{\bf Proposition.} {\it The commutation relations} \r{LL} {\it  for the Gauss
coordinates of the $L$-operator ${\cal L}^+(u,\eta)$}
$$
{\cal L}^+(u,\xi)=\left(
\begin{array}{cc}
h_+(u)/2&f_+(u)\\ e_+(u)&-h_+(u)/2
\end{array} \right)
$$
{\it are isomorphic to the commutation relations} \r{2.3}
{\it and
the generating functions $e_+(u)$, $f_+(u)$, $h_+(u)$
satisfy these commutation relations
if $\he_\la$, $\hf_\la$, $\hh_\la$ satisfy the
relations} \r{9}.
\smallskip

In order to prove this proposition we should use
the Fourier transform calculations and fix in  \r{LL}
either $\Im u_1<\Im u_2$ or $\Im u_1>\Im u_2$.

We conclude that the algebra $\ael$ is a classical limit of the algebra
 $\Ael$. Simple cal\-cu\-la\-ti\-ons
show that for $c =0$ the cobracket on $\ael$
 defined by the usual prescription
 \beq
\delta(x)=\lim{\tih\to 0}{\Delta(x)-\Delta'(x)\over\tih}
\label{Dd}
\eeq
coincides with cobracket \r{6}.

{}For $c\not =0$ we have no rights to define a cobracket $\delta(x)$ 
 by the prescription \r{Dd} since we cannot identify the tensor components
 in the image of $\Delta$. Nevertheless, formal application of the rule
 \r{Dd} yields the following result:
\bea
\delta(e_+(u))&=&  h_+(u)\wedge e_+(u)+
c\wedge\left({1\over2i}\,{de_+(u)\over du}-\eta^2{de_+(u)\over d\eta}
\right)\ ,\nn\\
\delta(f_+(u))&=&  f_+(u)\wedge h_+(u)+
c\wedge\left({1\over2i}\,{df_+(u)\over du}-\eta^2{df_+(u)\over d\eta}
\right)\ ,\nn\\
\delta(h_+(u))&=&2e_+(u)\wedge f_+(u)+
c\wedge\left({1\over2i}\,{dh_+(u)\over du}-\eta^2{dh_+(u)\over d\eta}
\right)\ .\label{2.4}
\eea
which can be also rewritten in  Fourier harmonics  as follows:
\bea
\delta \he_\la&=&-{1\over2}\vpint d\tau\ \hh_\tau\wedge\he_{\la-\tau}
\left(\cth\,\tau/2\eta+\th\,(\la-\tau)/2\eta \right)
+{\la\over2}\ \th\left({\la\over2\eta}\right) \he_\la\wedge c\ ,\nn\\
\delta \hf_\la&=&-{1\over2}\vpint d\tau\ \hf_{\la-\tau}\wedge\hh_\tau
\left(\cth\,\tau/2\eta+\th\,(\la-\tau)/2\eta \right)
+{\la\over2}\ \th\left({\la\over2\eta}\right) \hf_\la\wedge c\ ,\nn\\
\delta \hh_\la&=&-\intt d\tau\ \he_\tau\wedge\hf_{\la-\tau}
\left(\th\,\tau/2\eta+\th\,(\la-\tau)/2\eta \right)
+{\la\over2}\ \cth\left({\la\over2\eta}\right) \hh_\la\wedge c
 \ .\label{333}
\eea

\setcounter{equation}{0}
\section{The Rational Degeneration}

The affine algebras with a bialgebra structure usually appear as classical
doubles and thus are factorized into a sum of isotropic subalgebras.
It is not true for the algebra $\aelo$. Nevertheless,
as it  follows from the definition of the elements of the algebra
$\ael$ \r{11} and the generating functions \r{8a},
each substrip of the strips $\Pi^\pm$ defines a subalgebra of the
algebra $\ael$. It is clear from \r{11} that
in terms of Fourier components these subalgebras are
distinguished  by different asymptotics of the functions $g(\la)$,
$g'(\la)$, $g''(\la)$  at $\la\to\pm\infty$.
Let us consider the substrips $\tilde{\Pi}^\pm\subset \Pi^\pm$
\beq
\tilde{\Pi}^+=\left\{-{1\over2\eta}<\Im u<0\right\}
,\quad
\tilde{\Pi}^-=\left\{0<\Im u<{1\over2\eta}\right\}
\label{strips}
\eeq
and restrict the generating functions  $e^+(u)$,  $f^+(u)$,  $h^+(u)$
onto $\tilde{\Pi}^+$ and
$e^-(u)$,  $f^-(u)$,  $h^-(u)$ onto $\tilde{\Pi}^-$ respectively.
Let us denote corresponding subalgebras of $\ael$
as $\aelp$ and $\aelm$.
\bigskip

\unitlength 1mm
\linethickness{0.4pt}
\begin{picture}(140.00,60.00)
\put(139.67,10.00){\vector(1,0){0.2}}
\put(40.00,10.00){\line(1,0){99.67}}
\bezier{120}(40.33,11.00)(70.00,11.00)(70.00,11.00)
\bezier{96}(70.00,11.00)(86.00,12.33)(90.00,19.00)
\bezier{96}(90.00,19.00)(94.67,26.33)(110.00,27.00)
\bezier{120}(110.00,27.00)(123.67,28.33)(140.00,28.00)
\put(90.00,10.00){\line(0,1){5.00}}
\put(90.00,17.67){\line(0,1){4.33}}
\put(90.00,26.33){\line(0,1){3.8}}
\put(90.00,30.00){\line(1,0){5.00}}
\put(100.00,30.00){\line(0,1){0.00}}
\put(100.00,30.00){\line(1,0){5.00}}
\put(110.00,30.00){\line(1,0){5.00}}
\put(120.00,30.00){\line(1,0){5.00}}
\put(130.00,30.00){\line(1,0){7.33}}
\put(139.67,40.00){\vector(1,0){0.2}}
\put(140.00,40.00){\line(-1,0){99.67}}
\bezier{120}(139.67,41.00)(110.00,41.00)(110.00,41.00)
\bezier{96}(110.00,41.00)(94.00,42.33)(90.00,49.00)
\bezier{96}(90.00,49.00)(85.33,56.33)(70.00,57.00)
\bezier{120}(70.00,57.00)(56.33,58.33)(40.00,58.00)
\put(90.00,40.00){\line(0,1){5.00}}
\put(90.00,47.67){\line(0,1){4.33}}
\put(90.00,56.33){\line(0,1){3.8}}
\put(90.00,60.00){\line(-1,0){3.00}}
\put(80.00,60.00){\line(0,1){0.00}}
\put(80.00,60.00){\line(-1,0){5.00}}
\put(70.00,60.00){\line(-1,0){5.00}}
\put(60.00,60.00){\line(-1,0){5.00}}
\put(50.00,60.00){\line(-1,0){7.33}}
\put(16.33,20.00){\makebox(0,0)[cc]{\mbox{${1\over1+\ee^{-\la/\eta}}$}}}
\put(16.33,50.00){\makebox(0,0)[cc]{\mbox{${1\over1+\ee^{\la/\eta}}$}}}
\put(16.33,5.00){\makebox(0,0)[cc]{Fig.~2.}}
\end{picture}

Then in the limit $\eta\to+0$ these subalgebras become isotropic subalgebras
of the loop algebra $\widehat{\frak{sl}_2}$,
the generating
functions of these subalgebras
will be defined in  lower and upper half-planes and their
expressions via formal generators turns into the Laplace transform
(see Fig.~2.).

In this limit the family $\Ael$ turns
into the central extended Yangian double $\DY$ \cite{K,Ia} defined
in contrast to \cite{K,Ia} by means of $L$-operators $L^\pm(u)$
which are analytical in lower and upper half-planes respectively
and generated
by the continuous family  of the formal generators
$\he_\la$,  $\hf_\la$, $\hh_\la$
which are formal Fourier harmonics of the elements of the
$L$-operators $L^\pm(u)$ \cite{KLP2}. The Lie algebra $\ael$
turns at this limit into the central extension of $sl_2$-valued
 rational functions vanishing at $\infty$;
 subalgebras $\aelp$ and $\aelm$ turn into
subalgebras of rational functions analytical in lower and upper half-planes
respectively.

The elements of the algebra $\frak{a}_0(\widehat{\frak{sl}}_2)$
can be identified with $\frak{sl}_2$-valued distributions
\beq
e\ot {P(z)\over Q(z)},\quad
f\ot {P'(z)\over Q'(z)},\quad
h\ot {P''(z)\over Q''(z)}\ ,
\label{***}
\eeq
where $P(z)$, etc. are polynomials such that
${\rm deg}\,P(z)<{\rm deg}\,Q(z)$, ${\rm deg}\,P'(z)<{\rm deg}\,Q'(z)$ and
${\rm deg}\,P''(z)<{\rm deg}\,Q''(z)$.
Subalgebras
$\frak{a}^+_0(\widehat{\frak{sl}}_2)$ and
$\frak{a}^-_0(\widehat{\frak{sl}}_2)$   are distinguished by the
contours in  pairing the dis\-tri\-bu\-tions \r{***}
and the elements of the space of the basic functions on the complex plane
analytical in the lower and the upper half-planes respectively
and vanishing at the infinity. For the elements of the subalgebra
$\frak{a}^+_0(\widehat{\frak{sl}}_2)$ this contour is parallel to the
real axis and
lies above all zeros
of the polynomials $Q(z)$, $Q'(z)$ and $Q''(z)$ and vice versa for the
subalgebra $\frak{a}^-_0(\widehat{\frak{sl}}_2)$.

The central extention of the algebra $\frak{a}^+_0(\widehat{\frak{sl}}_2)$
is defined by the standard cocycle and the commutation relations between
the generating functions of the elements of the algebra
$\frak{a}^+_{0,c}(\widehat{\frak{sl}}_2)$ and can be formally obtained
from the commutation relations \r{3.1} at the limit $\eta\to+0$: 
\bea
{[}e_\pm(u_1), f_\pm(u_2){]}&=& {i\over u_1-u_2 }
\left(h_\pm(u_1)-h_\pm(u_2)\right)\ ,\nn\\
{[}h_\pm(u_1),h_\pm(u_2){]}&=& 0\ ,\nn\\
{[}h_\pm(u_1), e_\pm(u_2){]}&=& {2i\over u_1-u_2}\left(e_\pm(u_2)-
e_\pm(u_1)\right)\ ,\nn\\
{[}h_\pm(u_1), f_\pm(u_2){]}&=&{-2i\over u_1-u_2}\left(f_\pm(u_2)-
f_\pm(u_1)\right)\ ,\nn\\
{[}e_\pm(u_1), f_\mp(u_2){]}&=& {i\over u_1-u_2 }
\left(h_\mp(u_1)-h_\pm(u_2)\right)\mp{1\over (u_1-u_2)^2}\cdot c\ ,\nn\\
{[}h_\pm(u_1),h_\mp(u_2){]}&=& \mp{2\over(u_1-u_2)^2}\cdot c\ ,\nn\\
{[}h_\pm(u_1), e_\mp(u_2){]}&=& {2i\over u_1-u_2}\left(e_\mp(u_2)-
e_\pm(u_1)\right)\ ,\nn\\
{[}h_\pm(u_1), f_\mp(u_2){]}&=&{-2i\over u_1-u_2}\left(f_\mp(u_2)-
f_\pm(u_1)\right)\ .\nn
\eea
Subalgebras $\frak{a}^\pm_0(\widehat{\frak{sl}}_2)$ become isotropic
 and there is a nontrivial pairing between these subalgebras:
$$
\langle e_+(u_1),f_-(u_2)\rangle={i\over u_1-u_2},\quad
\langle h_+(u_1),h_-(u_2)\rangle={2i\over u_1-u_2}.
$$
The algebra $\frak{a}_0(\widehat{\frak{sl}}_2)$ is a bialgebra with
the cobracket:
\bea
\delta(e_\pm(u))&=&\pm h_\pm (u)\wedge e_\pm (u)\pm\frac{1}{2i}c\wedge
{d e_\pm (u)\over du}
\ ,\nn\\
\delta(f_\pm(u))&=&\pm f_\pm (u)\wedge h_\pm (u)\pm\frac{1}{2i}c\wedge {d f_\pm (u)\over du}
\ ,\nn\\
\delta(h_\pm(u))&=&\pm 2e_\pm (u)\wedge f_\pm (u)\pm\frac{1}{2i}c\wedge{d h_\pm (u)\over du}
\ .\nn
\eea

The representation \r{8a} in the form of the Fourier integrals
comes to be the Laplace transform (this can be visualized in Fig.1)
$$
e_\pm(u)=\pm\int_{0}^\infty d\la\ \theta(\la)\ee^{\mp i\la u}\ \he_{\mp\la},
\
f_\pm(u)=\pm\int_{0}^\infty d\la\ \theta(\la)\ee^{\mp i\la u}\ \hf_{\mp\la},
\
h_\pm(u)=\pm \int_{0}^\infty d\la\ \theta(\la)\ee^{\mp i\la u}\ \hh_{\mp\la},
$$
so the formal generators $\he_\la$,  $\hf_\la$, $\hh_\la$ at $\la\leq0$
form the subalgebra $\frak{a}^+_0(\widehat{\frak{sl}}_2)$
and $\he_\la$,  $\hf_\la$, $\hh_\la$ at $\la\geq 0$
the subalgebra $\frak{a}^-_0(\widehat{\frak{sl}}_2)$. The map
 \r{333} becomes a cobracket now
\bea
\delta \he_\la&=&\int_0^\la d\tau\ \bigl(\theta(\tau-\la)-\theta(\tau)\bigr)
\hh_\tau\wedge\he_{\la-\tau}
+\frac{\la}{2}{\rm sgn}(\la)\he_\la \wedge c
\ ,\nn\\
\delta \hf_\la&=&\int_0^\la d\tau\ \bigl(\theta(\tau-\la)-\theta(\tau)\bigr)
\hf_{\la-\tau}\wedge\hh_{\tau}
+\frac{\la}{2}{\rm sgn}(\la)\hf_\la \wedge c
\ ,\nn\\
\delta \hh_\la&=&2\int_0^\la d\tau\ \bigl(\theta(\tau-\la)-\theta(\tau)\bigr)
\he_{\tau}\wedge\hf_{\la-\tau}
+\frac{\la}{2}{\rm sgn}(\la)\hh_\la \wedge c
\ .\nn
\eea 
Here  $\theta(\la)$ is the step function
$$
\theta(\la)=\left\{\begin{array}{c}
1,\quad \la>0\\ 
1/2, \quad \la=0\\
0,\quad \la<0
\end{array}\right.
$$
and defines a Lie bialgebra structure on
$\frak{a}_0(\widehat{\frak{sl}}_2)$
at arbitrary central element $c$. This bialgebra is a classical double
of one of subalgebras
$\frak{a}^\pm_0(\widehat{\frak{sl}}_2)$.

\setcounter{equation}{0}
\section{The  Lie Algebra $\apssect$}

The structure of the Lie bialgebra $\ael$ (at level 0) from 
section 1 and section 2
can be automatically generalized to the elliptic case.
Using the Lie algebra of double periodic automorphic functions which
take values in $\frak{sl}_2$ \cite{RS} one can construct
on the half-parallelogram
of the periods $\Pi$ $(0,2K,iK',2K+iK')$
the Lie algebra of $\frak{sl}_2$-valued generalized functions
$\apo$ ($p=\ee^{\pi i\tau}=\ee^{-\pi K'/K}$).
Let us introduce the generating functions
$\sigma^{\pm}_a(u)$
of the elements of the
algebra $\apo$:
\beq
\sigma^+_a(u)=\sigma_a\ot\omega_a(z-u),
\quad  \sigma_a^-(u)=\sigma_a^+(u-K'),\quad a=1,2,3\ ,
\label{6.0}
\eeq
where $\sigma_a$ are the standard Pauli matrices
$$
\sigma_1=\left(\begin{array}{cc}0&1\\ 1&0\end{array}\right),\quad
\sigma_2=\left(\begin{array}{cc}0&-i\\ i&0\end{array}\right),\quad
\sigma_3=\left(\begin{array}{cc}1&0\\ 0&-1\end{array}\right)
$$
and functions $\omega_a(u)$ are ratios of elliptic Jacobi functions
of modulus $k$ \cite{BE}:
$$
\omega_1(u)={1\over{\rm sn}\,(u,k)},\quad
\omega_2(u)={{\rm dn}\,(u,k)\over{\rm sn}\,(u,k)},\quad
\omega_3(u)={{\rm cn}\,(u,k)\over{\rm sn}\,(u,k)}\ .
$$
Because of the addition theorem for elliptic functions
\beq
\omega_a(u-v)\omega_c(v)-
\omega_b(u-v)\omega_c(u)=
\omega_a(u)\omega_b(v)
\label{addthe}
\eeq
the generating functions $\sigma^\pm_a(u)$ satisfy the commutation
relations:
\bea
{[}\sigma_a^\pm(u_1),\sigma_b^\pm(u_2){]}&=&2\left[
i\omega_a(u_1-u_2)\sigma_c^\pm(u_2)-
i\omega_b(u_1-u_2)\sigma_c^\pm(u_1)\right]\ ,\nn\\
{[}\sigma_a^\pm(u_1),\sigma_b^\mp(u_2){]}&=&2\left[
i\omega_a(u_1-u_2)\sigma_c^\mp(u_2)-
i\omega_b(u_1-u_2)\sigma_c^\pm(u_1)\right]\ ,\label{4.1}
\eea
where $a,b,c$ are cyclic permutations of 1,2,3.
The cobracket
\beq
 \delta\sigma_a(u)=\sigma_b(u)\wedge\sigma_c(u)
\label{4.2}
\eeq
defines on $\apo$ a Lie bialgebra structure.

In analogy with the trigonometric case the algebra $\apo$ admits a
central extension given by the two-cocycle:
$$
B(x\otimes \varphi(z),y\otimes \psi(z))={1\over2K}
\int_{\partial \Pi}{dz\over2\pi i} \left({d\psi(z)\over d\tau}\varphi(z)-
\psi(z){d\varphi(z)\over d\tau}
\right) \ \langle x,y\rangle\ ,
$$
where $\partial \Pi$ is the boundary of the
period's half-parallelogram $(0,2K,iK',2K+iK')$. The  cocycle property
is a consequence of the addition theorem \r{addthe}. The values of
the cocycle on the generating functions are
\beq
B(\sigma_a^\pm(u_1),\sigma_b^\pm(u_2))={\delta_{a,b}\over K}
{d\omega_a(u_1-u_2)\over d\tau}\ ,
\label{BB}
\eeq
 the values of the cocycle for the other combinations of generating
 functions one can find applying to \r{BB} the automorphic conditions
 \r{6.0}. 
The central extended  Lie algebra $\ap$ is defined in terms of generating
functions subjected to the commutation relations \r{4.1} together with
 the relations \r{62} instead of trivial ones.
\beq
[\sigma^\pm_a(u_1),\sigma^\pm_a(u_2)]=
c\cdot
B(\sigma^\pm_a(u_1),\sigma^\pm_a(u_2)),\quad
[\sigma_a^\pm(u_1),\sigma^\mp_a(u_2)]=
c\cdot
B(\sigma_a^\pm(u_1),\sigma^\mp_a(u_2))\ .
\label{62}
\eeq

It should be  mentioned that
in terms of Fourier
harmonics
\bea
\sigma_1^+(u)&=&{i\pi\over K}\sum_{\ell\ {\rm odd}}\sigma_1^\ell
{p^\ell\over p^\ell-1}\ee^{2\pi i\ell u},\nn\\
\sigma_2^+(u)&=&{i\pi\over K}\sum_{\ell\ {\rm odd}}\sigma_2^\ell
{p^\ell\over p^\ell+1}\ee^{2\pi i\ell u},\nn\\
\sigma_3^+(u)&=&{i\pi\over K}\sum_{\ell\ {\rm even}}\sigma_3^\ell
{p^\ell\over p^\ell+1}\ee^{2\pi i\ell u},\nn
\eea
the commutation relations \r{4.1}, \r{62} become the relations for the
generators of central extended loop algebra $\widehat{sl_2}$:
\beq
[\sigma_a^k,\sigma_b^\ell]=2i\varepsilon_{abc}\sigma_c^{k+\ell}+
c\cdot B(\sigma_a^k,\sigma_b^\ell)
\label{sigma}
\eeq
with the standard cocycle
$$
B(\sigma_a^k,\sigma_b^\ell)=k\delta_{a,b}\delta_{k,-\ell}.
$$
But, as well as in the trigonometric case,
it is not sufficient to write down the relations \r{sigma} for the
description of the Lie algebra $\ap$. One should specify the elements of $\ap$
to be the series of the formal generators $\sigma_a^k$ with coefficients
of the type \r{11} and \r{12}, which appear after Fourier
decomposition of the generalized functions from the space
$\ap$. In terms of Fourier harmonics the cobracket \r{4.2}
is given by the relations
\bea
\delta(\sigma_1^k) &=& \sum_{i,j\atop i+j=k}
{p^k-1\over (p^i+1)(p^j+1)} \sigma_2^i\wedge\sigma_3^j\ ,\nn\\
\delta(\sigma_2^k) &=& \sum_{i,j\atop i+j=k}
{p^k+1\over (p^i+1)(p^j-1)} \sigma_3^i\wedge\sigma_1^j\ ,\nn\\
\delta(\sigma_3^k) &=& \sum_{i,j\atop i+j=k}
{p^k+1\over (p^i-1)(p^j+1)} \sigma_1^i\wedge\sigma_2^j\ .\nn
\eea

The central extended algebra $\ap$ can be written in $r$-matrix
formalism. Let 
\beq
r(u)=\sum_{a=1}^3 \omega_a(u)\sigma_a \ot\sigma_a
\label{r-ellip}
\eeq
be an elliptic solution of the classical Yang-Baxter equation
and
$$
{\cal L^+}(u)=\left(\begin{array}{cc}
\sigma^+_3(u) & \sigma^+_1(u)-i\sigma^+_2(u) \\
\sigma^+_1(u)+i\sigma^+_2(u) & -\sigma^+_3(u)
\end{array}\right),
$$
where $\sigma^+_a(u)$ satisfy the commutation relations \r{4.1}
and \r{62}. 
The $r$-matrix \r{r-ellip}
is traceless and the commutation relations \r{4.1} and \r{62} can
be written in the following form:
\beq
{[}{\cal L}^+_1(u_1),{\cal L}^+_2(u_2){]}={[}{\cal L}^+_1(u_1)+
{\cal L}^+_2(u_2),
r(u_1-u_2){]} +{1\over K}
{d r(u_1-u_2)\over d\tau}\cdot c\ .
\label{565}
\eeq
This representation allows us to demonstrate that the Lie algebra
$\ap$ is a classical limit of the quantum elliptic algebra $\Apq$
\cite{MJ}. Indeed, in the `$RLL$' formalism the algebra $\Apq$  
is defined by the
relations
\beq
R(v;\tih,\tilde\tau)L^+_1(u_1)L^+_2(u_2)=L^+_2(u_2)L^+_1(u_1)
R(v;\tih,\tilde\tau^*)\ ,
\label{63}
\eeq
where $\tih$ is a deformation parameter and $R(v;\tih,\tilde\tau)$
is the Baxter elliptic $R$-matrix \cite{B}
$$
\rho(v)\left(\begin{array}{cccc}
a(v)&0&0&d(v)\\ 0&b(v)&c(v)&0\\ 0&c(v)&b(v)&0\\  d(v)&0&0&a(v)
\end{array}\right).
$$
The normalization factor
$\rho(v)$ provides
the crossing symmetry and unitarity conditions and
\bea
a(v;\tih,\tilde\tau)&=&\sn(\tih+iv),\quad
b(v;\tih,\tilde\tau)\ =\ \sn(iv),\nn\\
c(v;\tih,\tilde\tau)&=&\sn(\tih),\quad
d(v;\tih,\tilde\tau)\ =\ \tilde k \,\sn(\tih)\sn(iv)\sn(\tih+iv).
\label{68}
\eea
In \r{68} $\sn(x)=\sn(x,\tilde k)$ is
Jacobi's elliptic function of modulus $\tilde k$. The modular parameter
$\tilde \tau=i {\tilde K}'/\tilde K$
is defined by the half-periods $\tilde K$ and ${\tilde K}'$ and
related to $\tilde k$ in the standard way. In \r{63} $\tilde\tau^*$
is the following expression:
\beq
\tilde\tau^*=\tilde\tau+\alpha\tih c\ ,
\label{64}
\eeq
where $c$ is the central element and $\CC$-number $\alpha$ depends on
the half-period $\tilde K$ and modulus $\tilde k$.

To obtain the commutation relations \r{565} from \r{63} first we
have to use 
the Landen transform 
$$
k'={1-\tilde k\over1+\tilde k},\quad 
\zeta={\tih(1+\tilde k)\over 2},\quad u={i(1+\tilde k)v}
$$
applied to the matrix elements of the $R$-matrix
and then the imaginary Jacobi
transform which connect the elliptic functions of argument $iu$ and
supplementary
modulus $k'$ with those of argument $u$ and 
modulus $k$ (see details in \cite{B}).
The $R$-matrix in this new parametrization reads:
$$
R(u;\zeta,\tau)=\tilde\rho(u)\left[
1+\sum_{a=1}^3 W_a(u)\sigma_a\ot\sigma_a\right],
$$
where
$$
W_1={\sn(\zeta,k)\over\sn(u+\zeta,k)},\quad
W_2={\sn(\zeta,k)\,\dn(u+\zeta,k)\over\sn(u+\zeta,k)\,\dn(\zeta,k)},\quad
W_3={\sn(\zeta,k)\,\cn(u+\zeta,k)\over\sn(u+\zeta,k)\,\cn(\zeta,k)}
$$
and had been used by E.~Sklyanin in \cite{Sk}.
Since Landen transform do not affect much the modular parameter $\tau$
the prescription of the central extension \r{64} will be the same:
$\tau^*=\tau+\zeta c/K$.

After all the transformations
the classical limit of the quantum relations \r{63} means that
$\zeta\to0$. In this limit $L^+(u)=1+\zeta\,{\cal L}^+(u)+o(\zeta^2)$,
$R(u;\zeta,\tau)=\left[1+\zeta r(u,\tau) + \zeta^2
r_1(u,\tau) + o(\zeta^3)\right]$,
$R(u;\zeta,\tau^*)=\left[1+\zeta r(u,\tau) + \zeta^2
\tilde r_1(u,\tau) + o(\zeta^3)\right]$
 and 
$$
\tilde r_1(u,\tau)-r_1(u,\tau)={c\over K}{dr(u,\tau)\over d\tau}
$$
so \r{565} appears as first nontrivial coefficient of the Taylor expansion
of \r{63} with respect to $\zeta$. Since classical $r$-matrix is
traceless the scalar factor $\tilde\rho(u)$ does not contribute
into calculation of \r{565} (compare with subsection 4.2).

\section{Acknowledgement}

The authors would like to acknowledge the discussions with A.~Gerasimov,
A.~Laptev and
A.~Rosly.  
S.~Pakuliak would like to thank Prof.~T.~Miwa for the invitation 
to visit RIMS.
The work was supported in part by grants RFBR-96-01-01421
(V.~Tolstoy),
RFBR-96-01-00814a (S.~Kho\-rosh\-kin), 
RFBR-96-01-01101 (D.~Lebedev),
RFBR-96-02-19085 (S.~Pakuliak),
INTAS-93-0166-Ext (D.~Lebedev),
INTAS-93-2058-Ext  (S.~Pakuliak)
and by Award  No. RM2-150 of the U.S. Civilian Research \& Development
Foundation (CRDF) for the Independent States of the Former Soviet Union
(S.~Khoroshkin, D.~Lebedev, S.~Pakuliak).

\bigskip

\noindent
{\Large\bf Note added in the replaced version}
\bigskip

\noindent
 One can notice that a map $\delta$ given by the relation \r{333},
 respects the structure of Lie algebra $\ael$ but
 does not satisfy (co)Jacoby identity and thus furnish the algebra
$\ael$ with a structure of quasi Lie bialgebra. Moreover,
 cobracket \r{333} can be also defined via classical $r$-matrix, described
 in \cite{semenov}, \cite{jimbo2}:
\bea
\delta(x_\la)&=&{[}x_\la\ot1+1\ot x_\la, r],\quad x=\he,\hf,\hh,\nn\\
r&=&
\intt d\la\ {\hf_{-\la}\ot\he_{\la}+\he_{-\la}\ot\hf_{\la}\over
1+\ee^{\la/\eta}}
+{1\over2}  \vpint d\la\ {\hh_{-\la}\ot\hh_{\la}\over
1-\ee^{\la/\eta}} + {1\over2}(c\ot d + d\ot c),\nn
\eea
where operator $d$ is a gradation operator 
$$
[d, x_\la]=\la x_\la\quad\mbox{for} \quad x=\he,\hf,\hh.
$$
This is in a correspondence with results \cite{jimbo1}, \cite{jimbo2}
 establishing an equivalence of elliptic affine algebra with a twisted
 quantum affine algebra.

\renewcommand{\refname}{}

\end{document}